\documentclass[final]{siamltex}
\usepackage{graphicx, amssymb, amsmath}
\usepackage{color}
\usepackage{overpic,subeqnar}
\usepackage{soul} 
\usepackage{booktabs}
\usepackage{boxedminipage}
\usepackage{amsfonts}
\usepackage{algorithmic}
\usepackage{algorithm,booktabs}
\usepackage[mathscr]{eucal}

\begin{document}
\title{Online interpolation point refinement for reduced order models using a genetic algorithm \thanks{J. N. Kutz acknowledges support from the Air Force Office of Scientific Research (FA9550-15-1-0385).}}

%\title{Online genetic algorithm selection of nearly optimal interpolation points in parametrized, nonlinear model order reduction \thanks{J. N. Kutz acknowledges support from the Air Force Office of Scientific Research (FA9550-15-1-0385).}}

\author{
Syuzanna Sargsyan\footnote{University of Washington, Department of Applied Mathematics, Seattle, WA 98195, USA
\texttt{ssuzie@uw.edu}}, 
Steven L. Brunton\footnote{University of Washington, Department of Mechanical Engineering, Seattle, WA 98195, USA
\texttt{sbrunton@uw.edu}}, 
J. Nathan Kutz \footnote{University of Washington, Department of Applied Mathematics, Seattle, WA 98195, USA
\texttt{kutz@uw.edu}}
%
%
%Dante Kalise\footnote{Johann Radon Institute for Computational and Applied Mathematics, 
%Altenbergerstra\ss e 69, 4040 Linz, Austria. 
%\texttt{dante.kalise@oeaw.ac.at}} 
}

\maketitle

\begin{abstract}
A genetic algorithm procedure is demonstrated that refines the selection of interpolation points of the
discrete empirical interpolation method (DEIM) when used for constructing reduced order models for time dependent and/or parametrized nonlinear partial differential equations (PDEs) with proper orthogonal decomposition.  The method achieves nearly optimal interpolation
points with only a few generations of the search, making it potentially useful for {\em online} refinement of the sparse
sampling used to construct a projection of the nonlinear terms.  
With the genetic algorithm, points are optimized to jointly minimize reconstruction error and enable dynamic regime classification.  
The efficiency of the method is demonstrated on 
two canonical nonlinear PDEs:  the cubic-quintic Ginzburg-Landau equation and the Navier-Stokes equation for flow around a cylinder.  Using
the former model, the procedure can be compared to the ground-truth optimal interpolation points, showing that 
the genetic algorithm quickly achieves nearly optimal performance and reduced the reconstruction error by nearly an order of
magnitude.  
\end{abstract}



\maketitle

%%%%%%%%%%%%%%%%%%%%%%%%%%%%%%%%%%%%%%%%%%%%%%%%%%%%%
\section{Introduction}
\label{Section1}
\setcounter{section}{1}
\setcounter{equation}{0}
\setcounter{theorem}{0}
\setcounter{algorithm}{0}
\renewcommand{\theequation}{\arabic{section}.\arabic{equation}}

Scientific computation has become critically enabling in almost every field of scientific and engineering study, enabling
simulations with modern computers that were thought to be out of reach even a decade ago and suggesting the possibility of exascale computing architectures.  The continued scaling of memory, processor speed and parallelization enable studies of increasingly sophisticated multi-scale physical systems.  Despite these advances, significant challenges and computational bottlenecks still remain in efficiently computing dynamics of extremely high-dimensional systems, such as high Reynolds turbulent flow.  Reduced order models (ROMs) are of growing importance as a critically enabling mathematical framework for reducing the dimension of such large systems.
The core of the ROM architecture relies on two key innovations:  (i) The POD-Galerkin method~\cite{HLBR_turb,karen1}, which is used for projecting the high-dimensional nonlinear dynamics to a low-dimensional subspace in a principled way, and (ii) sparse sampling (gappy POD) of the state space for interpolating the nonlinear terms required for the subspace projection.   
The focus of this manuscript is on a sparse sampling innovation for ROMs. Specifically, a method for optimizing sampling locations for both reconstruction and identification of parametrized systems.   We propose an algorithm comprised of two components:  (i) an {\em offline} stage that produces initial sparse sensor locations, and (ii) an {\em online} stage that uses a short, genetic search algorithm for 
producing nearly optimal sensor locations.  The technique optimizes for both reconstruction error and classification efficacy, leading to an attractive {\em online} modification of commonly used gappy POD methods.
 
The importance of sparse sampling of high-dimensional systems, especially those manifesting low-dimensional dynamics, 
was recognized early on in the ROMs community.  Thus sparse sampling has already been established as a critically enabling mathematical framework for model reduction through methods such as gappy POD and its variants~\cite{gap1,gap2,karni,Carlberg:2013}.  
More generally, sparsity promoting methods are of growing importance in physical modeling and scientific computing~\cite{Wang2011prl,Schaeffer2013pnas,Ozolicnvs2013pnas,mackey2014compressive,Bai2014aiaa,Brunton2016pnas}.  
 The seminal work of Everson and Sirovich~\cite{gap1} first established how the gappy POD could play a transformative role in the mathematical sciences.  In their sparse sampling scheme, random measurements were used to perform reconstruction tasks of inner products.  Principled selection of the interpolation points, through the gappy POD infrastructure~\cite{gap1,gap2,karni,Carlberg:2013} or missing point (best points) estimation (MPE)~\cite{patera,mpe}, were quickly incorporated into ROMs to improve performance.  More recently, the transformative empirical interpolation method (EIM)~\cite{eim} and its most successful variant, the POD-tailored discrete empirical interpolation method (DEIM)~\cite{deim}, have provided a greedy, sparse samplimg algorithm that allows for nearly optimal reconstructions of nonlinear terms of the original high-dimensional system.  The DEIM approach combines projection with interpolation.  Specifically, the DEIM uses selected interpolation indices  to specify  an interpolation-based projection for a nearly optimal $\ell_2$ subspace approximating the nonlinearity.  

It is well-known that the various sparse sampling techniques proposed are not optimal, but have been shown to be sufficiently robust to provide accurate reconstructions of the high-dimensional system.  The DEIM algorithm has been particularly successful
for nonlinear model reduction of time-dependent problems.  Interestingly, for parametrized systems, the DEIM algorithm needs
to be executed in the various dynamical regimes considered, leading to a library learning mathematical framework~\cite{sargsyan2015}.
Thus efficient sparse sampling locations for both classification and reconstruction can be computed in an {\em offline} manner across various dynamical regimes.  Again, they are not optimal, but they are robust for building ROMs.  We build upon the DEIM library learning framework~\cite{sargsyan2015}, showing that nearly-optimal sparse sampling can be achieved with a short {\em online}, genetic algorithm search from the learned DEIM libraries.  This improves both the classification and reconstruction accuracy of the sparse sampling, making it an attractive performance enhancer for ROMs.  

The paper is outlined as follows:  In Sec.~\ref{sec:ROM}, the basic ROM architecture is outlined.
Section~3 reviews the various innovations of the sparse sampling architecture, including the library building
procedure used here.  Section~4 develops the genetic search algorithm for {\em online} improvement of
sparse sampling locations.  The method advocated here is demonstrated in Sec.~5 on two example problems:
the complex cubic-quintic Ginzburg-Landau equation and fluid flow past a circular cylinder.  Concluding remarks are
provided in Sec.~6. 

\section{Reduced Order Modeling}
\label{sec:ROM}

In our analysis, we consider a parametrized, high-dimensional system of nonlinear differential equations that arises, for example, from the finite-difference discretization of a partial differential equation.  In the formulation proposed, the linear and nonlinear terms
for the state vector ${\bf u}(t)$ are separated:
 \begin{equation}
 \label{eq:complex}
   \frac{d {\bf u}(t)}{dt}=L{\bf u}(t)+N({\bf u}(t),\mu),
 \end{equation}
 where ${\bf u}(t)=[u_1(t) \,\, u_2(t) \,\, \cdots \,\, u_n(t)]^T\in \mathbb{R}^n$ and $n\gg 1$.  Typically, $u_j (t) = u(x_j, t)$ is the value of the field of interest discretized at the spatial location $x_j$.
 The linear part of the dynamics is given by the linear operator $L\in \mathbb{R}^{n\times n}$ and the nonlinear terms are
 in the vector $N({\bf u}(t)) = [N_1({\bf u}(t),\mu) \quad N_2({\bf u}(t),\mu) \quad \cdots \quad N_n({\bf u}(t),\mu)]^T\in \mathbb{R}^n$. 
 The nonlinear function is evaluated component-wise at the $n$ spatial grid points used for 
 discretization.  Note that we have assumed, without loss of generality, that the parametric dependence $\mu$ is in the nonlinear term. 
 
Typical discretization schemes for achieving a prescribed spatial resolution and accuracy require the number of discretization 
points $n$ to be very large, resulting in a high-dimensional state vector.  For sufficiently complicated problems where significant spatial refinement
is required and/or higher spatial dimension problems (2D or 3D computations, for instance) can potentially lead to a
computationally intractable problem where ROMs are necessary.    The POD-Galerkin method~\cite{HLBR_turb,karen1} 
is a principled dimensionality-reduction scheme that approximates the function ${\bf u}(t)$ with 
rank-$r$-optimal basis functions where $r\ll n$.  
These optimal basis functions are computed from a singular value decomposition of a time series
of snapshots of the {{nonlinear dynamical}} system (\ref{eq:complex}).   Given the snapshots of the
state variable ${\bf u}(t_j)$ at $j=1, 2, \cdots , p$ times, the snapshot matrix ${\bf X}=[{\bf u} (t_1) \,\, {\bf u}(t_2) \,\, \cdots \,\, {\bf u}(t_p)]
  \in \mathbb{R}^{n\times p}$ is constructed and the SVD of ${\bf X}$ is computed
\begin{equation}
  {\bf X}={\bf \Psi \Sigma W}^* \, .
\end{equation}
 The $r$-dimensional basis for optimally approximating 
 ${\bf u}(t)$ is given by the first $r$ columns of matrix ${\bf \Psi}$, denoted by ${\bf \Psi}_r$.
 The POD-Galerkin approximation is then computed from the following decomposition
 \begin{equation} {\bf u}(t)\approx {\bf \Psi}_r {\bf a}(t)
 \label{eq:podG}
 \end{equation} 
where ${\bf a}(t)\in \mathbb{R}^r$ is the time-dependent coefficient vector and $r\ll n$. 
Note that such an expansion implicitly assumes a separation of variables between time and space.
 Plugging this modal expansion into the governing equation (\ref{eq:complex}) and applying orthogonality (multiplying
 by ${\bf \Psi}_r^T$) gives the dimensionally reduced evolution
   \begin{equation}
   \label{eq:pod}
     \frac{d{\bf a}(t)}{dt}={\bf \Psi}_r^T L {\bf \Psi}_r {\bf a}(t)+ {\bf \Psi}_r^T N({\bf \Psi}_r {\bf a}(t),\mu)
   \end{equation}
where  it is noted that ${\bf \Psi}_r^T {\bf \Psi}_r = {\bf I}$, the $r\times r$ identity matrix.

By solving this system of much smaller dimension, the solution of a high-dimensional {{nonlinear dynamical}}
system can be approximated with optimal (in an $\ell_2$ sense) basis functions.  
Of critical importance is evaluating the nonlinear terms in an efficient way using
the gappy POD or DEIM mathematical architecture.  Otherwise, the evaluation of the nonlinear terms still requires calculation of functions and inner products with the original dimension $n$.  In certain cases, such as for the quadratic nonlinearity in the Navier-Stokes equations, the nonlinear terms can be computed once in an offline manner.  However, parametrized systems generally require repeated evaluation of the nonlinear terms as the POD modes change with $\mu$.  Equation (\ref{eq:pod}) provides the starting point of the
ROM architecture.  It also illustrates one of the two critical innovations of ROMs:  a principled dimensionality reduction.

\section{Gappy Sampling and Discrete Empirical Interpolation Method}

In and of itself, the dimensionality reduction provided by the POD method does not ensure that 
the reduced model (\ref{eq:pod}) remains low-dimensional.  This is simply due to the evaluation of
the nonlinear term ${\bf \Psi}_r^T N({\bf \Psi}_r {\bf a}(t),\mu)$.   Specifically, the nonlinearity forces
the evaluation of the nonlinear function and inner products at every time step when advancing (\ref{eq:pod}) with
a time-stepping scheme.  Such inner products require computations that scale with $n$, keeping the ROM fixed
in the higher dimensional space.  In contrast, the linear term ${\bf \Psi}_r^T L {\bf \Psi}_r$ can be computed
once, either at the beginning of the simulation, or in an offline stage.  This product of matrices yields
a matrix of size $r\times r$, thus requiring computations that scale with $r$ to update the linear portion of the evolution.

To avoid the costly computations associated with approximating $\mathbf{N}=N({\bf \Phi}_r {\bf a}(t))$, we compute
an approximation to the nonlinearity through projection and interpolation instead of evaluating it directly.  
A considerable reduction in complexity is achieved by the DEIM, for instance, because evaluating the approximate nonlinear term does not require a prolongation of the reduced state variables back to the original high dimensional state approximation required to evaluate the nonlinearity in the POD approximation.  The DEIM therefore improves the efficiency of the POD approximation and achieves a complexity reduction of the nonlinear term with a complexity proportional to the number of reduced variables.   The DEIM constructs these specially selected interpolation indices that specify an interpolation-based projection to provide a nearly $\ell_2$ optimal subspace approximation to the nonlinear term without the expense of orthogonal projection~\cite{deim}.

\subsection{Sparse (Gappy) Sampling}

To better understand how the computation of the nonlinearity through projection and interpolation, we consider
taking $m$ measurements ($m>r$ with $O(m)\sim O(r)$) of the full state vector ${\bf u}$.
Specifically, only $m\ll n$ measurements are required for reconstruction, allowing us to 
define the sparse representation variable $\tilde{\bf u}\in\mathbb{R}^{m}$
\begin{equation}
  \tilde{\bf u} = {\bf P} {\bf u}
  \label{eq:gap}
\end{equation}
where the measurement matrix ${\bf P}\in\mathbb{R}^{m\times n}$ specifies $m$ measurement locations 
of the full state ${\bf u}\in\mathbb{R}^{n}$.
As an example, the measurement matrix might take the form
\begin{equation}
\label{eq:P}
  {\bf P} = \left[  \begin{array}{cccccccccc}  
     1 & 0          &  \cdots     &     &      & & & & \cdots  & 0 \\
     0 & \cdots   & 0             &  1 &  0  & \cdots & & & \cdots  & 0 \\
     0 & \cdots   &                &     & \cdots     & 0 & 1 & 0 &\cdots & 0 \\
     \vdots &  0   &                &     & \cdots     & 0 & 0 & 1 &\cdots & \vdots \\
     0 & \cdots   &                &     & \cdots     & 0 & 0 & 0 &\cdots & 1 
  \end{array}    \right]
\end{equation}
where measurement locations take on the value of unity and the matrix elements are zero elsewhere.
The matrix ${\bf P}$ defines a projection onto an $m$-dimensional measurement space $\tilde{\bf u}$ that is
used to approximate ${\bf u}$.  

The sparse sampling of~(\ref{eq:gap}) forms the basis of the {\em Gappy POD} method introduced
by Everson and Sirovich~\cite{gap1}.   In their example of eigenface reconstruction, they used a
small number of measurements, or gappy data, to reconstruct an entire face.  This serves as the
basis for approximating the nonlinear inner products in (\ref{eq:pod}) and overcoming the computational
complexity of the POD reduction.  The measurement matrix ${\bf P}$ allows for an approximation of the state vector ${\bf u}$
from $m$ measurements.  The approximation is given by substituting (\ref{eq:podG}) into (\ref{eq:gap}):
\begin{equation}
  \tilde{\bf u} \approx {\bf P} \sum_{j=1}^m \tilde{ a}_j {\psi}_j 
  \label{eq:gapapprox}
\end{equation}
where the coefficients $\tilde{a}_j$ minimize the error in approximation
 \begin{equation}
   \| \tilde{\bf u}-{\bf P u}\|.
   \label{eq:error}
 \end{equation}
The goal is to determine the $\tilde{a}_j$ despite the fact that taking inner products
of (\ref{eq:gapapprox}) can no longer be performed.  Specifically, the vector $\tilde{\bf u}$ is of
length $m$ whereas the POD modes are of length $n$.
Implied in this sparse sampling is that the modes $\psi_j(x)$ are in general not orthogonal 
over the $m$-dimensional support of $\tilde{\bf u}$, which is denoted as $s[\tilde{\bf u}]$. 
Specifially, the following two relationships hold
\begin{subeqnarray}
&& M_{jk} = \left( \psi_j,\psi_k \right)  = \delta_{jk}  \\
&& M_{jk} = \left( \psi_j,\psi_k \right)_{s[\tilde{\bf u}]}  \neq 0 \,\,\,\, \mbox{for all} \, j, k
\label{eq:M}
\end{subeqnarray}
where $M_{jk}$ are the entries of the Hermitian matrix ${\bf M}$ and $\delta_{jk}$ is the 
Kroenecker delta function.  The fact that the POD modes
are not orthogonal on the support $s[\tilde{\bf u}]$ leads us to consider alternatives
for evaluating the vector $\tilde{\bf a}$.

As demonstrated by Everson and Sirovich~\cite{gap1}, the $\tilde{a}_j$ is determined by  
a least-square fit error algorithm.  Specifically, we project the full state vector ${\bf u}$ onto the support space and
determine the vector $\tilde{\bf a}$ with the equation:
\begin{equation}
  {\bf M} \tilde{\bf a} = {\bf f} 
\end{equation}
where the elements of the matrix ${\bf M}$ are given by (\ref{eq:M}b) and the components of the vector $f_k$ are given by
$  f_j =  \left( {\bf u},  \psi_j \right)_{s[\tilde{\bf u}]} $.  The pseudo-inverse for determining $\tilde{\bf a}$ is a least-square fit
algorithm.
Note that in the event the measurement space is sufficiently dense, or as the support space is the entire space,
then ${\bf M}={\bf I}$ and $\tilde{\bf a} \rightarrow {\bf a}$, thus implying the eigenvalues of ${\bf M}$ approach unity as the number of measurements become dense.  Once the vector $\tilde{\bf a}$ is determined, then a reconstruction of the solution
can be performed using
\begin{equation}
  {\bf u}(x,t) \approx  {\bf \Psi} \tilde{\bf a} \, .
  \label{eq:separate2}
\end{equation}
It only remains to consider the efficacy of the measurement matrix ${\bf P}$.  Originally, random measurements were proposed~\cite{gap1}.  However, the ROMs community quickly developed principled techniques based upon, for example, minimization of the
condition number of ${\bf M}$~\cite{gap2}, selection of minima or maxima of POD modes~\cite{karni}, and/or greedy algorithms of EIM/DEIM~\cite{eim,deim}.  Thus $m$ measurement locations were judiciously chosen for the task of accurately interpolating the nonlinear
terms in the ROM.  This type of sparsity has been commonly used throughout the ROMs community.

\subsection{DEIM Algorithm}

The DEIM algorithm constructs two low-rank spaces through the SVD:  one for the full system using the snapshots matrix ${\bf X}$, and
a second using the snapshot matrix composed of samples of the nonlinearity alone.  Thus a low-rank representation of the
nonlinearity is given by
\begin{equation}
  {\bf N} = \left[ N({\bf u}_1) \,\, N({\bf u}_2) \,\, \cdots \,\, N({\bf u}_p)  \right] = {\bf \Xi \Sigma}_N {\bf W}_N^*
\end{equation}
where the matrix ${\bf \Xi}$ contains the optimal (in an $\ell_2$ sense) basis set for spanning the nonlinearity.
Specifically, we consider the rank-$m$ basis set  ${\bf \Xi}_m=[\boldsymbol{\xi}_1 \,\, {\boldsymbol \xi}_2 \,\, \cdots \,\, {\boldsymbol \xi}_m]$
that approximates the nonlinear function ($m\ll n$ and $m\sim r$).  
The approximation to the nonlinearity $\mathbf{N}$ in this SVD basis set is given by:
\begin{equation}
 \mathbf{N}\approx {\bf \Xi}_m {\bf c}(t)
\end{equation}
where ${\bf c}(t)$ is similar to ${\bf a}(t)$ in (\ref{eq:podG}).
Since this is a highly overdetermined system, a suitable vector ${\bf c}(t)$ can be found by selecting only $m$ 
rows of the system.  The DEIM algorithm provides a greedy search algorithm for selecting an
appropriate $m$ rows.  Although not guaranteed to be optimal, in practice the row selection tends to provide
sampling points that are accurate for reconstruction of the full state.
 
The DEIM algorithm begins by considering the vectors ${\bf e}_{\gamma_j}\in\mathbf{R}^n$ which are the $\gamma_j$-th column of 
the $n$ dimensional identity matrix.  We can then construct the projection matrix 
${\bf P}=[{\bf e}_{\gamma_1} \,\, {\bf  e}_{\gamma_2} \,\, \cdots \,\, {\bf e}_{\gamma_m}]$ which is chosen so that ${\bf P}^T{\bf \Xi}_m$ is nonsingular. 
Then ${\bf c}(t)$ is uniquely defined from 
\begin{equation}
{\bf P}^T\mathbf{N} ={\bf P}^T {\bf \Xi}_m {\bf c}(t),
\end{equation}
and thus, 
\begin{equation}
\label{eq:nonlinearity}
\mathbf{N}\approx {\bf \Xi}_m ({\bf P}^T {\bf \Xi}_m)^{-1} {\bf P}^T\mathbf{N}. 
\end{equation}
The tremendous advantage of this result for nonlinear model reduction is that the term ${\bf P}^T\mathbf{N}$   
 requires evaluation of nonlinearity only at $m$ indices, where $m\ll n$. 
The DEIM further proposes a principled method for choosing the basis vectors ${\boldsymbol \xi}_j$ and 
indices $\gamma_j$.  The DEIM algorithm, which is based upon a greedy-like search, is detailed in~\cite{deim} and
further demonstrated in Table~\ref{table:alg}.

\begin{table*}[t]
   \caption{ \label{table:alg} DEIM algorithm for finding approximation basis for the nonlinearity and interpolation indices.}
\begin{center}
\begin{tabular}[c]{|p{6cm}|p{6cm}|}
 \hline
  \multicolumn{2}{ |c| } { \bf  DEIM algorithm}\\  \hline
  \multicolumn{2}{ |c| } {\bf  {   Basis}} \\  \hline
  $\bullet$ construct snapshot matrix & $\mathbf{X}=[\mathbf{u}(t_1) \,\, \mathbf{u}(t_2) \,\, \cdots \,\, \mathbf{u}(t_p)]$\\  \hline
  $\bullet$ construct nonlinear snapshot matrix & $\mathbf{N}=[N({\bf u}(t_1)) \,\, N({\bf u}(t_2)) \,\, \cdots \,\, N({\bf u}(t_p))]$\\  \hline
  $\bullet$ singular value decomposition of $\mathbf{N}$ & $\mathbf{N}={\bf \Xi}{\bf \Sigma}_N {\bf W}_N^*$\\  \hline
  $\bullet$ rank-$m$ approximating basis & ${\bf \Xi}_m=[{\boldsymbol \xi}_1 \,\, {\boldsymbol \xi}_2 \,\, \cdots \,\, {\boldsymbol \xi}_m]$\\ \hline
  \multicolumn{2}{ |c| } {\bf  {   Interpolation Indices (Iteration Loop)}} \\  \hline
  $\bullet$ choose the first index (initialization) %corresponds to the entry of $v_1$ that has the largest   magnitude 
  & $[\rho, \gamma_1]=\max|{\boldsymbol \xi}_1|$\\ \hline
  $\bullet$ approximate ${\boldsymbol \xi}_j$ by ${\boldsymbol \xi}_1,...,{\boldsymbol \xi}_{j-1}$  at indices $\gamma_1,...,\gamma_{j-1} $
  %$v_j$ is a linear combination of $v_1, ...v_{j-1}$ 
  &  Solve for ${\bf c}$:
  ${\bf P}^T {\boldsymbol \xi}_j={\bf P}^T{\bf \Xi}_{j-1}{\bf c}$ with ${\bf P}=[{\bf e}_{\gamma_1} \,\, \cdots \,\, {\bf e}_{\gamma_{j-1}}]$  \\ \hline
  $\bullet$ select $\gamma_j$ and loop ($j=2, 3, ..., m$) & $[\rho,\gamma_j ]=\max|{\boldsymbol \xi}_j-{\bf \Xi}_{j-1} {\bf c}|$\\ \hline    
  \end{tabular}   
\end{center} 
 \end{table*}

\subsection{Application to Library Learning for parametrized ROMs}

The DEIM algorithm is {highly effective} for determining sampling (sensor) locations.  
Such sensors can be used with sparse representation and
compressive sensing to (i) identify dynamical regimes, (ii) reconstruct
the full state of the system, and (iii) provide an efficient nonlinear model reduction and POD-Galerkin prediction
for the future state.   Given the parametrized nature of the evolution equation (\ref{eq:pod}), we use the
concept of library building  which arises in machine learning from leveraging low-rank ``features" from data.  In the ROM community, it has recently become an issue of intense investigation.  Indeed, a variety of recent works~\cite{sargsyan2015,siads,Bright:2013,epj,amsallem1,amsallem2,karen2,karen3,karen4,Kaiser2014jfm} have produced libraries of ROM models that can be selected and/or interpolated  through measurement and classification.  Before these more formal
 techniques based upon machine learning were developed, it was already realized that parameter domains
 could be decomposed into subdomains and a local ROM/POD computed in each subdomain.
 Patera {\em et al.}~\cite{sub1} used a partitioning based on a binary tree whereas 
 Amsallem {\em et al.}~\cite{sub2} used  a Voronoi Tessellation of the domain.  Such methods
 were closely related to the work of Du and Gunzburger~\cite{sub3} where the data snapshots were 
 partitioned into subsets and multiple reduced bases computed.  Thus the concept of library building 
 is well established and intuitively appealing for parametrized systems.

We capitalize on these recent innovations and build optimal interpolation locations from multiple dynamics states~\cite{sargsyan2015}.  However, the focus of this work is on computing, in an online fashion, nearly optimal sparse sensor locations from interpolation points found to work across all the libraries in an offline stage.  The offline stage uses the DEIM architecture as this
method gives good starting points for the interpolation.  The genetic algorithm we propose then improves upon the interpolated points
by a quick search of nearby interpolation points.  It is the pre-computed library structure and interpolation points that
allow the genetic algorithm to work with only a short search.

\section{Genetic Algorithm for Improved Interpolation}

\begin{figure}[t]
\centering
\vspace*{-1in}
\hspace*{-.8in}
\begin{overpic}[width=1.3\textwidth]{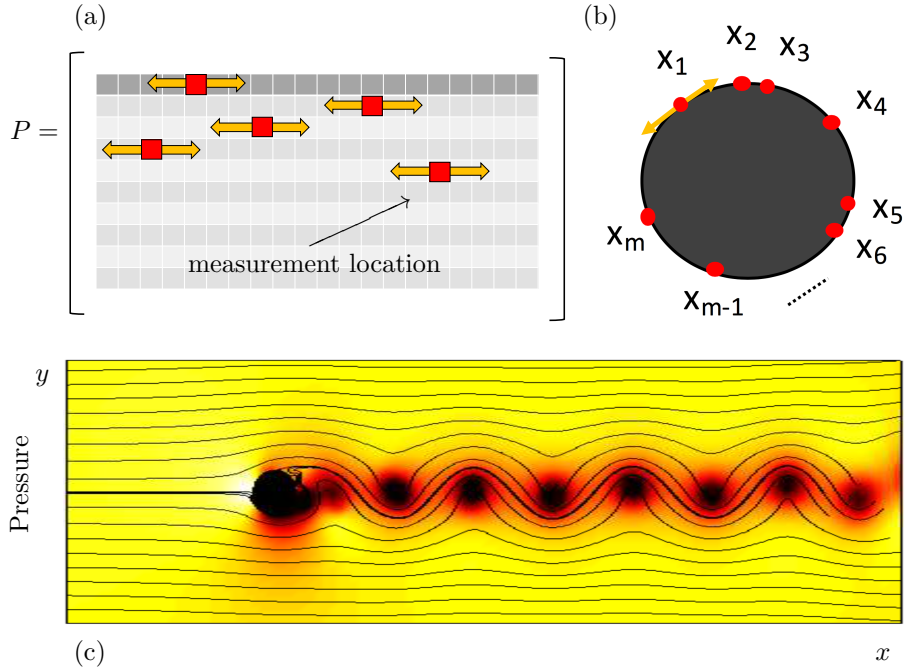}
\put(20,62){(a)}
\put(60,62){(b)}
\put(20,12){(c)}
\put(15,22){\rotatebox{90}{Pressure}}
\put(15,53){$P=$}
\put(17,34){$y$}
\put(83,12){$x$}
\put(38,44){\rotatebox{25}{$\xrightarrow{\hspace*{1.3cm}}$}}
\put(29,42.5){measurement location}
\end{overpic}
\vspace*{-1in}
\caption{Illustration of the genetic algorithm search for optimal sparse sampling locations.  (a) The measurement matrix
$P$ is constructed from the DEIM algorithm across libraries of POD modes~\cite{sargsyan2015}.  Given that the
matrix $P$ has already been demonstrated to produce good interpolation points,  the interpolation indices are
shifted to nearby locations in a genetic algorithm search for the best interpolation configuration.  (b)  The algorithm
has a practical physical interpretation for the example of flow around a cylinder where the sensor locations are shifted
for best performance.  (c)  The near-optimal interpolation points found from the genetic algorithm are used to 
both classify the dynamic regime, in this example it discovers the Reynolds number, and to reconstruct the full pressure field
(on and off the cylinder) with minimal residual~\cite{Bright:2013}. \label{fig:ga}}
\end{figure}

The background Secs.~2 and 3 provide the mathematical framework for the innovations of
this paper.  Up to this point, the DEIM architecture for parametrized PDEs~\cite{sargsyan2015} provides
good interpolation points for the ROM method.  Our goal is to make the interpolation points optimal
or nearly so.  Unfortunately, non-convex problems such as this are extremely difficult
to optimize, leading to the consideration of genetic algorithms, which are a subset of evolutionary algorithms,
for determining near optimal interpolation points.

The genetic algorithm principal is quite simple:  given a set
of feasible trial solutions (either constrained or unconstrained), an objective (fitness)
function is evaluated.  The idea
is to keep those solutions that give the minimal value of the objective
function and mutate them in order to try and do even better.  Mutation in our context involves
randomly shifting the locations of the interpolation points.  
Beneficial mutations that give a better minimization, such as good classification and
minimal reconstruction error, are
kept while those that perform poorly are discarded.
The process is repeated through a number of
iterations, or {\em generations}, with the idea that better and better fitness function values 
are generated via the mutation process.

More precisely, the genetic algorithm can be framed as the 
constrained optimization problem with the objective function
\begin{equation}
\min \|\tilde{\mathbf{u}} - \mathbf{P}\mathbf{u}\|_2 \quad \text{subject to correct classification}
%  \min f({\bf P}) 
\end{equation}
where ${\bf P}$ is a measurement matrix used for interpolation.  Suppose that $m$ mutations,
as illustrated in Fig.~\ref{fig:ga}, are given for the matrix ${\bf P}$ so that
\begin{equation}
 j^{\text{th}}\text{ guess is}\,\, {\bf P}_j  \, .
\end{equation}
Thus $m$ solutions are evaluated and compared with each other in
order to see which of the solutions generate the smallest objective
function since our goal is to minimize it.  We can order the guesses
so that the first $p<m$ gives the smallest values of $f({\bf P})$.  
Arranging our data, we then have
\begin{eqnarray}
  &  \mbox{keep}    & {\bf P}_j \,\, j=1, 2, \cdots, p \\
  &  \mbox{discard}\,\,\,\, & {\bf P}_j \,\, j=p+1, p+2 , \cdots , m \nonumber \, .
\end{eqnarray}
Since the first $p$ solutions are the best, these are kept in
the next generation.  In addition,we now generate $m-p$ new trial 
solutions that are randomly mutated from the $p$ best solutions.
This process is repeated through a finite number of iterations $M$
with the hope that convergence to the optimal, or near-optimal, solution is achieved.
Table~\ref{table:alg2} shows the algorithm structure particular to our application.
In our simulations, $m=100$ mutations are produced and $p=10$ are kept for 
further mutation.  The number of generations is not fixed, but we find that even with $M=3$,
significant improvement in reconstruction error can be achieved.

As we will show, the algorithm provides an efficient and highly effective method for
optimizing the interpolation locations, even in a potentially online fashion.  A disadvantage
of the method is that there are no theorems guaranteeing that the iterations will converge
to the optimal solution, and there are many reasons the genetic search can fail.
However, we are using it here in a very specific fashion.   Specifically, our initial 
measurement matrix ${\bf P}$ is already quite good for classification and reconstruction
purposes~\cite{sargsyan2015}.  Thus the algorithm starts close to the optimal solution.
The goal is then to further refine the interpolation points so as to potentially cut down on
the reconstruction and classification error.  The limited scope of the algorithm mitigates
many of the standard pitfalls of the genetic algorithm.

\begin{table*}[t]
   \caption{ \label{table:alg2} DEIM algorithm for finding approximation basis for the nonlinearity and interpolation indices.}
\begin{center}
\begin{tabular}[c]{|p{8cm}|p{4cm}|}
 \hline
  \multicolumn{2}{ |c| } { \bf  Genetic search algorithm}\\  \hline
  $\bullet$ construct initial measurement matrix ~\cite{sargsyan2015} & $\bf P$\\  \hline
  $\bullet$ perturb measurements and classify &${\bf P} \rightarrow {\bf P}_1, {\bf P}_2, {\bf P}_3, \cdots, {\bf P}_m $ \\  \hline
  $\bullet$ keep matrices with correct classification & ${\bf P}_k, {\bf P}_j, {\bf P}_\ell, \cdots$ \\  \hline
  $\bullet$ save ten best measurement matrices & ${\bf P}\rightarrow {\bf P}_1, {\bf P}_2, \cdots, {\bf P}_{10}$ \\ \hline
  $\bullet$ repeat steps for $M$ generations & \\ \hline
  $\bullet$ randomly choose one of ten best ${\bf P}$ and repeat &  \\ \hline
  \end{tabular}   
\end{center} 
 \end{table*}

%Added by Susie:\\
%Here is how GA works:\\
%1. Start with location from DEIM PRE (or DEIM+1 PRE)\\
%2. Perturb the location and do classification (without noise in measurements)\\
%3. If classified correctly, compute the error and save the location \\
%4. Go to step 2 and do it 100 times\\
%5. Save the best 10 (smallest reconstruction error) locations from the above algorithm\\
%6. Chose number of generations $M$ \\
%7. 100 times uniformly at random choose one of the ten locations and repeat steps 2 and 3\\
%8.  Save the best ten locations and go to step 7 $M$ times.

\section{Model Problems}\label{sec:model}

Two models help illustrate the principles and success of the genetic search algorithm coupled with DEIM.  In the first
example, only three interpolation points are necessary for classification and reconstruction~\cite{sargsyan2015}.  Moreover,
for this problem, a brute force search optimization can be performed to rigorously identify the best possible interpolation points.
This allows us to compare the method advocated to a ground truth model.  In the second example, the classical problem
of flow around a cylinder is considered.  This model has been ubiquitous in the ROMs community for demonstrating 
new dimensionality-reduction techniques.

\subsection{Cubic-Quintic Ginzburg-Landau Equation}

The Ginzburg-Landau (GL) equation is on of the canonical models of applied mathematics and mathematical physics as it manifests a wide range dynamical behaviors~\cite{Cross:1993}.  It is widely used in the study of pattern forming systems, bifurcation theory and dynamical systems.  Its appeal stems from its widespread use in the sciences:  modeling phenomena as diverse as condensed matter physics to biological waves.  The particular variant considered here is the cubic-quintic GL with fourth-order diffusion~\cite{kutz:SIAM}:
  \begin{equation}
    iU_t+\left(\displaystyle \frac{1}{2}-i\tau\right)U_{xx}-i\kappa U_{xxxx}+\left(1-i\mu\right) |U|^2U+
(\nu-i\epsilon )|U|^4U-i\gamma U=0, 
\label{eq:gl}
  \end{equation}
% \begin{eqnarray}
% && \hspace{-.2in} i {U}_t + \left(  \frac{1}{2} - i\tau  \right) {U}_{xx} - i \kappa {U}_{xxxx}
% + (1- i\mu) |{U}|^2 {U}  \nonumber \\
% &&  \hspace{.2in}+ (\nu - i \varepsilon) |{U}|^4 {U} - i\gamma {U} \! =\! 0,
%  \label{eq:gl}
%\end{eqnarray}
%
where ${U}(x, t)$ is a complex valued function of space and time.   Under discretization
of the spatial variable, $U(x,t)$ becomes a vector ${\bf u}$ with $n$ components, i.e. 
${\bf u}_j (t)=U(x_j,t)$ with $j=1, 2, \cdots n$.

\begin{figure}[t]
\centering
\begin{overpic}[width=1.0\textwidth]{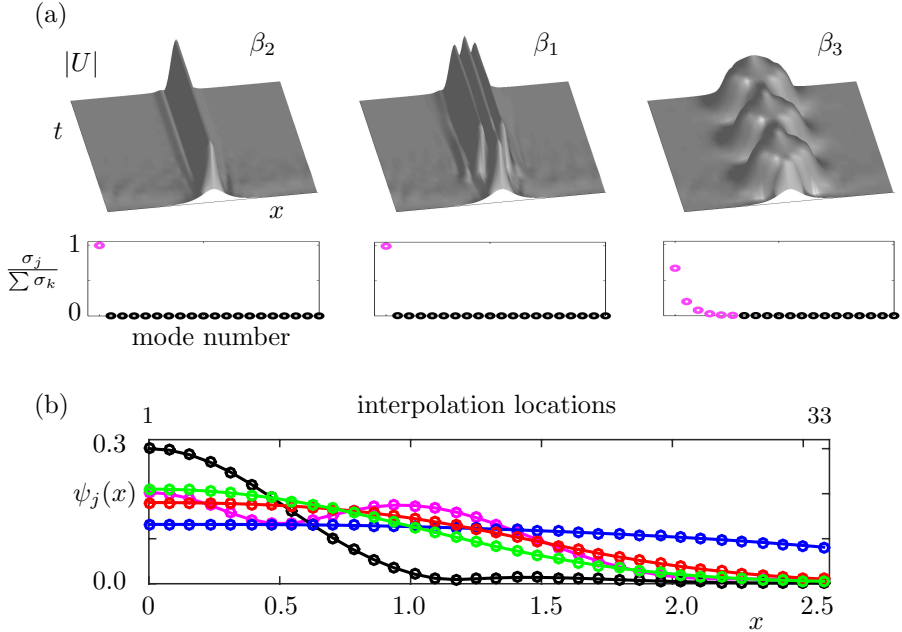}
\put(16,29){\small{1}}
\put(84,29){\small{33}}
\put(16,9.5){0 \hspace*{0.48in} 0.5 \hspace*{0.40in} 1.0 \hspace*{0.45in} 1.5 \hspace*{0.44in} 2.0 \hspace*{0.44in} 2.5}
\put(38,30){interpolation locations}
\put(5,70){(a)}
\put(5,30){(b)}
\put(8,65){$|U|$}
\put(29,50){$x$}
\put(7,58){$t$}
\put(27,67){$\beta_2$}
\put(56,67){$\beta_1$}
\put(85,67){$\beta_3$}
\put(2,44){$\frac{\sigma_j}{\sum{\sigma_k}}$}
\put(8.5,46.5){$1$}
\put(8.5,39.5){$0$}
\put(15,37){mode number}
\put(11,26){0.3}
\put(11,12){0.0}
\put(9,21){$\psi_j(x)$}
\put(78,8){$x$}
%\put(-4,-6){$x=0$}
%\put(-2,-3){\rotatebox{80}{$\xrightarrow{\hspace*{0.3cm}}$}}
%\put(13,-4){spatial grid point}
\end{overpic}
\vspace*{-0.50in}
\caption{(a) Evolution dynamics of (\ref{eq:gl}) for three different parameter ($\beta$) regimes as
highlighted in Table~\ref{table:vals}.  The intensity evolution is illustrate for $t\in[0,40]$ on the domain $x\in[-10,10]$ on a domain discretized with $n=1024$ points.  The subfigure for each evolution denotes the decay of singular values with the magenta
dots representing the modes retained in our library $\Psi_L$.  The modes retained account for 99.9\% of 
the total variance. (b) Profile of 5 POD modes $\psi_j(x)$ retained in the library of modes $\Psi_L$.   The figure shows only a small portion of the domain, $x\in[0,2.6]$, in order to highlight the 33 interpolation points that are possible to use over this selected domain.  Note that a majority of the modal structure is
contained in this domain, which suggests that this restricted domain may contain optimal interpolation points for the gappy POD evaluation of nonlinear terms for model reduction.  Indeed, these 33 interpolation points contain the interpolation points discovered by DEIM as well as the optimal points discovered
by exhaustive search (See Figs.~\ref{fig:gl1}-\ref{fig:gl3}).}
 \label{fig:cqg}
\end{figure}

An efficient and exponentially accurate numerical solution to (\ref{eq:gl}) can be found using standard spectral 
methods~\cite{Kutz:2013}.  Specifically, the equation is solved by Fourier transforming in the spatial
dimension and then time-stepping with an adaptive 4th-order Runge-Kutta method.  The extent of the spatial
domain is $x\in[-20,20]$ with $n=1024$ discretized points.   Importantly, in what follows the interpolation indices
are dictated by their position away from the center of the computational domain.  The center of the domain is 
at $x_{0}$ which is the 513th point in the domain.  The interpolation indices demonstrated are relative
to this center point.

\begin{table}[tb] 
\begin{center}
 \begin{tabular}[c]{|c|c|c|c|c|c|c|c|} %{|p{0.2cm}|p{0.4cm}|p{0.4cm}| p{0.4cm}| p{0.5cm}| p{0.4cm}|p{0.4cm}|p{0.4cm}| p{0.5cm}|p{2cm}| p{2cm}| p{2cm}| }

\hline
&$\tau$&$\kappa$ &$\mu$ & $\nu$ & $\epsilon$ & $\gamma$ & description\\
\hline
$\beta_1$ & -0.3    &-0.05     &1.45    &0     &-0.1  &-0.5& 3-hump, localized\\
\hline
$\beta_2$ & -0.3    &-0.05     &1.4  &0   &-0.1  &-0.5  &  localized, side lobes\\
\hline
$\beta_3$ & 0.08    & 0      &0.66      &-0.1    &-0.1     &-0.1    &breather\\
\hline
$\beta_4$ & 0.125  &0    &1    &-0.6    &-0.1    &-0.1    & exploding soliton\\
\hline
$\beta_5$ & 0.08    &-0.05    &0.6    &-0.1    &-0.1    &-0.1    & fat soliton\\
\hline
$\beta_6$ & 0.08    &-0.05    &0.5    &-0.1    &-0.1    &-0.1    & dissipative soliton\\
\hline
   
 \end{tabular}
 \end{center}
 \caption{Values of the parameters from equation (\ref{eq:gl}) that lead to six distinct dynamical regimes. To exemplify our
 algorithm, the first, third and fifth regimes
 will be discussed in this paper.}
 \label{table:vals}
\end{table}

To generate a variety of dynamical regimes, the parameters of the cubic-quintic GL are tuned to
a variety of unique dynamical regimes.  The unique parameter regime considered are denoted by
the parameter  $\beta=(\tau, \kappa, \mu, \nu, \epsilon, \gamma )$ which indicates the specific 
values chosen.  
Table \ref{table:vals} shows six different parameter regimes that have unique low-dimensional attractors as 
described in the table.  It has been shown in previous work that only three interpolation points are
necessary for classification of the dynamical state, state reconstruction and future state prediction~\cite{sargsyan2015,siads}.  
This previous work also explored how to construct the sampling matrix ${\bf P}$ from the DEIM algorithm
and its multiple dynamical state.

\begin{figure}[t]
\centering
\begin{overpic}[width=0.8\textwidth]{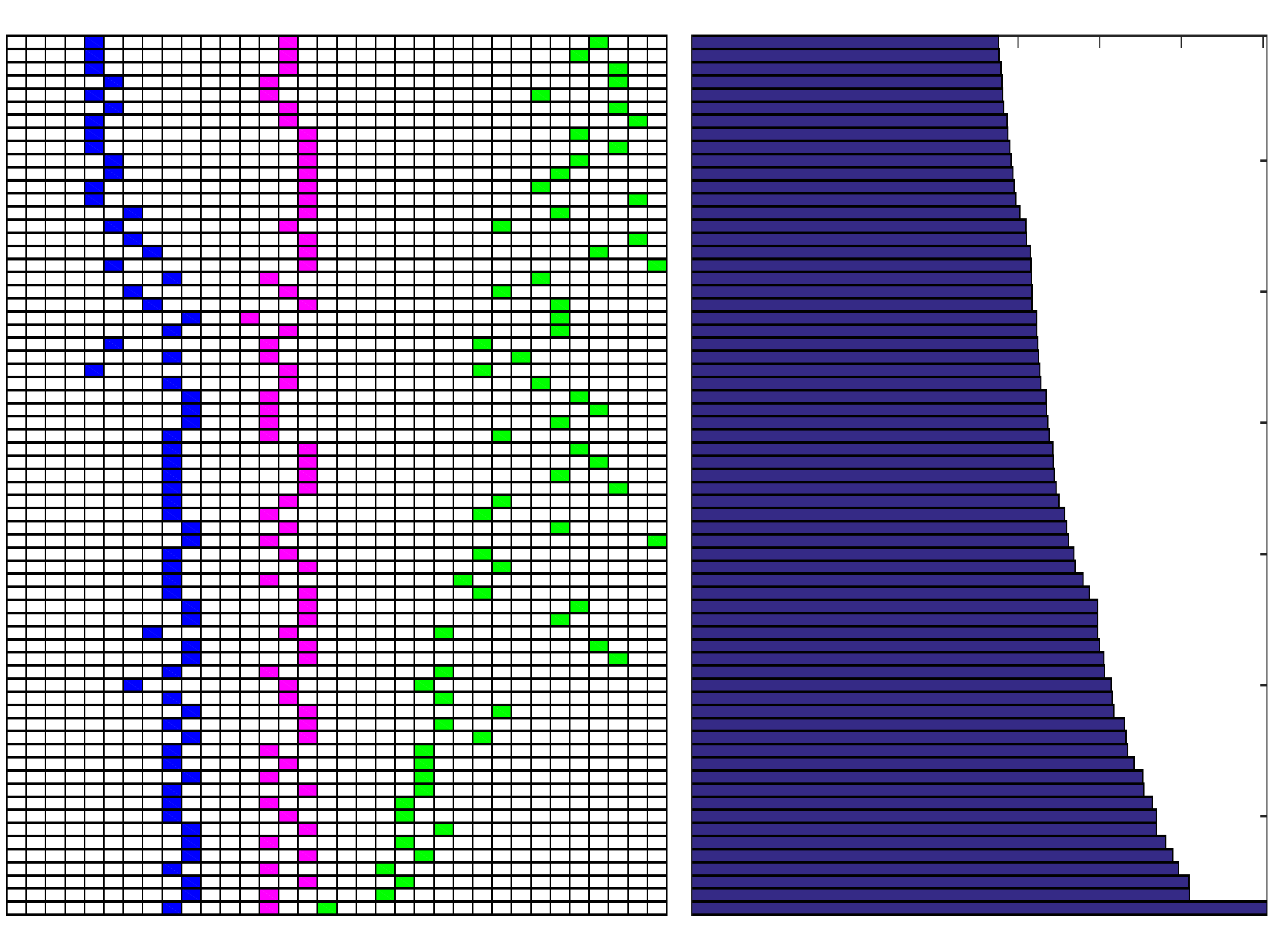}
\put(55,74){\small{0}}
\put(0.5,74){\small{1}}
\put(50,74){\small{33}}
\put(97,74){\small{0.2}}
\put(75,74){\small{0.1}}
\put(70,79){Error}
\put(9,79){interpolation locations}
\put(0,79){(a)}
\put(55,79){(b)}
\put(-4,-6){$x=0$}
\put(-2,-3){\rotatebox{80}{$\xrightarrow{\hspace*{0.3cm}}$}}
\put(13,-4){spatial grid point}
\put(-5,15){\rotatebox{90}{best interpolation triplets $\xrightarrow{\hspace*{0.5cm}}$}}
\end{overpic}
\vspace*{0.2in}
\caption{Results of an exhaustive brute force search for selecting the best three interpolation point triplets that correctly classify the dynamical regimes in the absence of noisy measurements.   From this subset, white noise is added to the measurements and 400 rounds of classification are performed.  Only the measurement triplets giving above 95\% accuracy for the classification of each dynamical regime are retained.  The retained triplets are then sorted by the reconstruction error as shown in (a).  The corresponding error is shown in (b).  Note that the interpolation indices are selected from
the first 33 points as shown in Fig.~\ref{fig:cqg}.}
 \label{fig:gl1}
\end{figure}
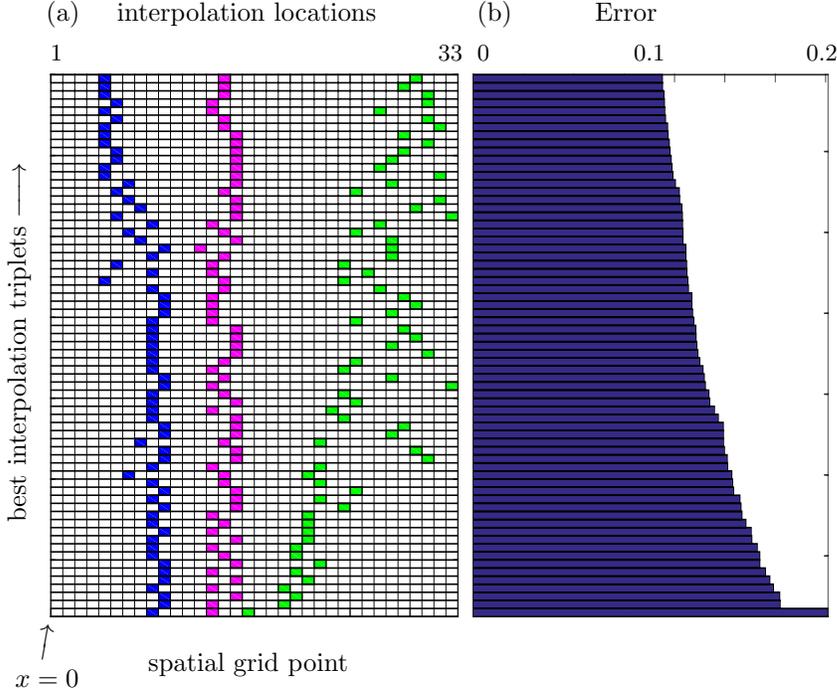

We will execute the genetic algorithm outlined in Table~\ref{table:alg2} for improving the sampling matrix ${\bf P}$ initially
determined from the algorithm in \cite{sargsyan2015}.  Before doing so, we consider a brute force search of the best possible
three measurement locations based upon their ability to classify the correct dynamical regime and minimize reconstruction error.  Although generally this is an $np$-hard problem, the limited number of sensors and 
small number of potential locations for sensors allow us to do an exhaustive search for the best interpolation locations.
The brute force search first selects all indices triplets (selected from interpolation points 0 to 33 as suggested by~\cite{sargsyan2015}) that correctly classify the dynamical regimes in the absence of noisy
measurements.   From this subset, white noise is added to the measurements and 400 rounds of classification 
are performed.  Only the measurement triplets giving above 95\% accuracy for the classification of each dynamical regime are retained.  The retained triplets are then sorted by the reconstruction error.  Figure~\ref{fig:gl1} shows the triplet
interpolation points retained from the exhaustive search with the classification criteria specified and the position of the interpolation points along with the reconstruction error.    The DEIM algorithm proposed in \cite{sargsyan2015} produces interpolation
points with reconstruction errors nearly an order of magnitude larger than those produced from the exhaustive search.
Our objective is to use a genetic algorithm to reduce our error by this order of magnitude and produce interpolation
points consistent with some of the best interpolation points displayed in Fig.~\ref{fig:gl1}.

\begin{figure}[t]
\hspace*{-.2in}
\begin{overpic}[width=0.55\textwidth]{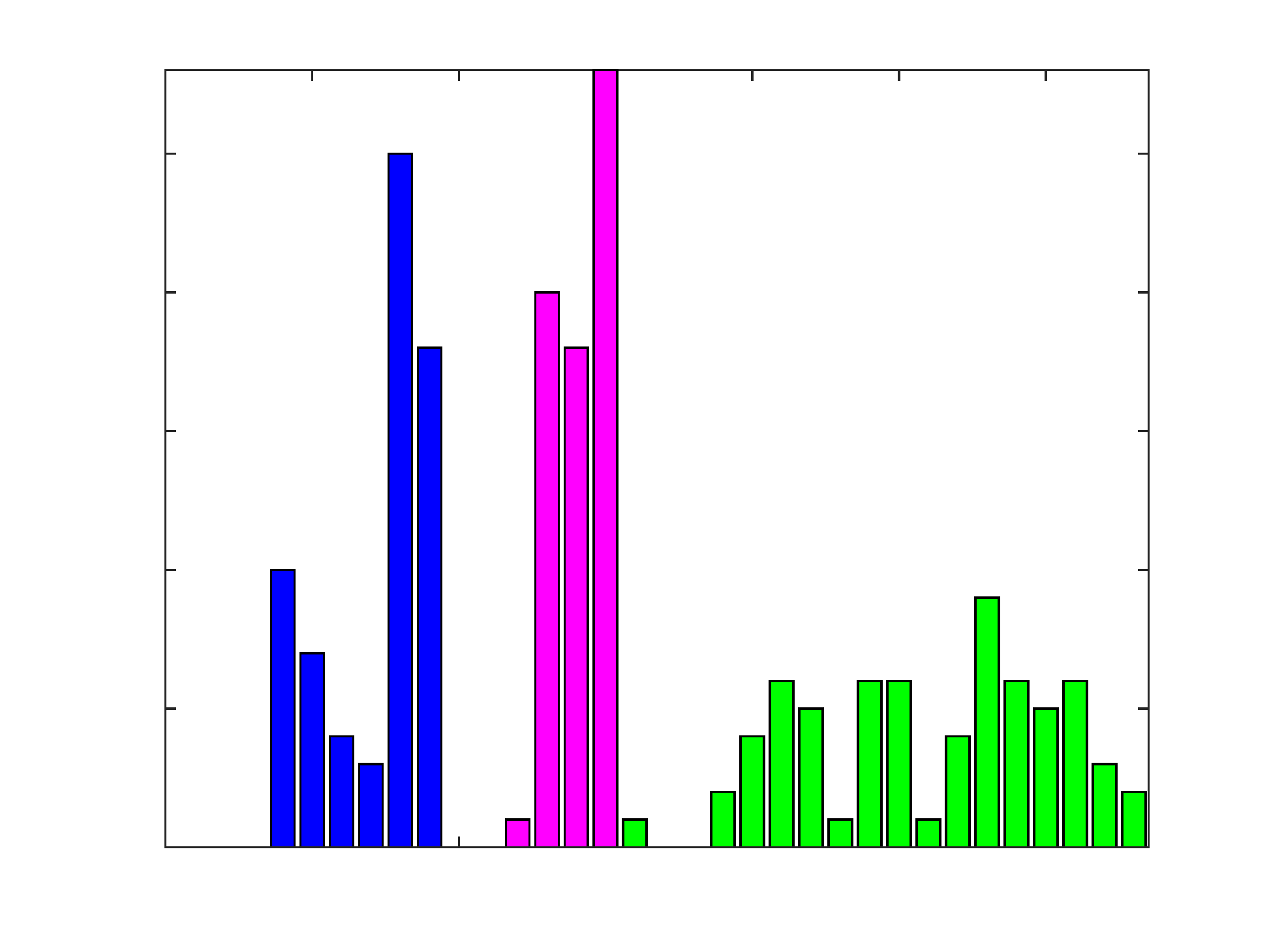}
\put(17,59){(a)}
\put(25,-1){interpolation locations}
\put(12,2){1}
\put(87,2){33}
\put(4,27){\rotatebox{90}{\# sensors}}
%\put(50,55){$\beta_3$}
\end{overpic}
\hspace*{-.3in}
\begin{overpic}[width=0.55\textwidth]{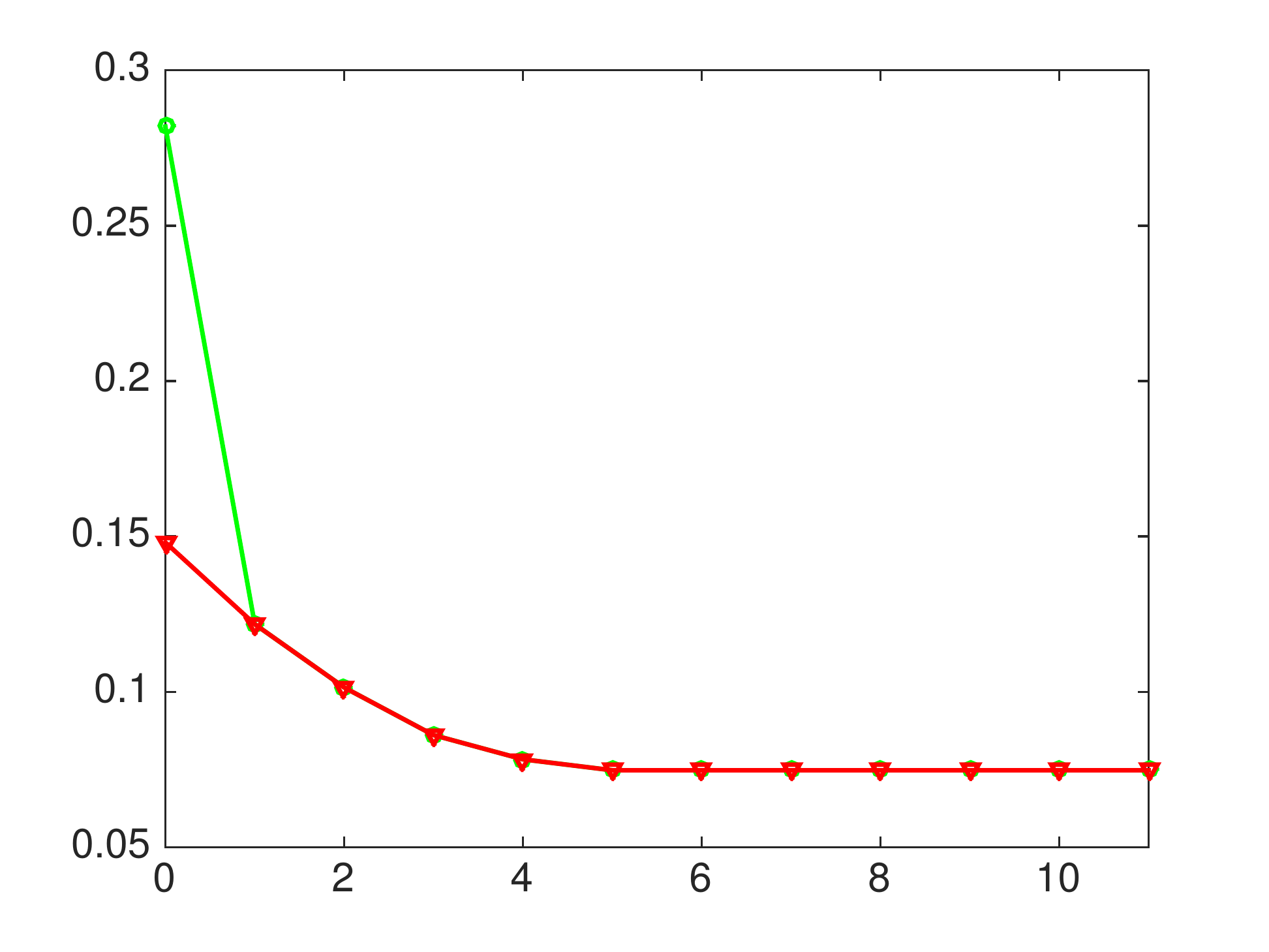}
\put(17,59){(b)}
\put(5,50){$E$}
\put(35,0){generation \#}
%\put(50,55){$\beta_3$}
\end{overpic}
\caption{(a) Histogram of the first (blue), second (magenta) and third (green) interpolation
points found from the exhaustive optimization search algorithm in Fig.~\ref{fig:gl1}.  
Note that the interpolation indices are selected from
the first 33 points as shown in Fig.~\ref{fig:cqg}.
The first two interpolation
points are highly localized near key spatial locations.  Note that the histograms demonstrate that
the first interpolation point selected from the standard DEIM algorithm is not an optimal point.  Instead,
for this case we can consider the 2nd-4th point instead, shifting the selection of optimal points by one iteration
and resulting in the DEIM+1 algorithm, i.e. we generate the first 4 interpolation points from DEIM and drop the first.  
(b) The genetic algorithm is executed starting form the sparse sampling matrix ${\bf P}$ of DEIM and DEIM+1 showing
that within 2-5 generations nearly optimal interpolation points are found in regards to the error $E$.
}
 \label{fig:gl2}
\end{figure}
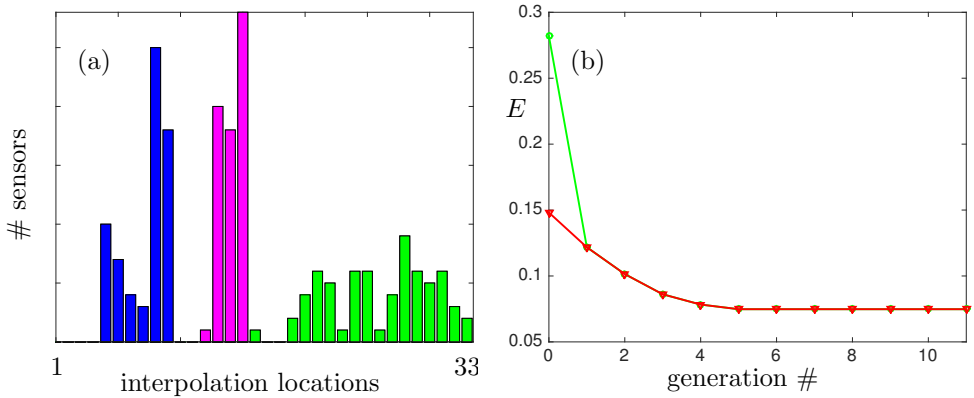

The brute for search produces a number of interpolation triplets whose reconstruction accuracy are quite similar.
Clearly displayed in the graph is the clustering of the interpolation points around critical spatial regions.  A histogram
of the first (blue), second (magenta) and third (green) interpolation points is shown in Fig.~\ref{fig:gl2}(a).
The first two interpolation points have a narrow distribution around the 4th-8th interpolation points and 12th-16th
interpolation points respectively.  The third interpolation point is more diffusely spread across a spatial region with
improvements demonstrated for interpolation points further to the right in Fig.~\ref{fig:gl1}(a).  This histogram provides critical information about sensor and interpolation point locations.  Of note is the fact that the DEIM algorithm
always picks the maximum of the first POD mode as an interpolation location.  This would correspond to a measurement
at $x=0$.  However, none of the candidate triplets retained from a brute force search consider this interpolation point
to be important. In fact, the interpolation points starting from the second iteration of the DEIM algorithm are what seem
to be important according to the brute force search.  This leads us to conjecture that we should initially use the triplet
pair from the 2nd-4th DEIM points rather than the 1st-3rd DEIM points.  We call these the DEIM+1 interpolation points as
we shift our measurement indices to the start after the first iteration of DEIM.

The genetic algorithm search can now be enacted from both the DEIM locations computed in ~\cite{sargsyan2015} and
the DEIM+1 locations suggested by the exhaustive search.   Figure~\ref{fig:gl2}(b) shows the convergence of the
genetic search optimization procedure starting from both these initial measurement matrices ${\bf P}$.  It should
be noted that the DEIM+1 initially begins with approximately half the error of the standard DEIM algorithm, suggesting
it should be used as a default.  In this specific scenario, both initial measurement matrices are modified and converge
to the near-optimal solution within only 3-5 generations of the search.  This is a promising result since the mutations
and generations are straightforward to compute and can potentially be done in an online fashion.  The benefit from
this approach is a reduction of the error by nearly an order of magnitude, making it an attractive scheme.

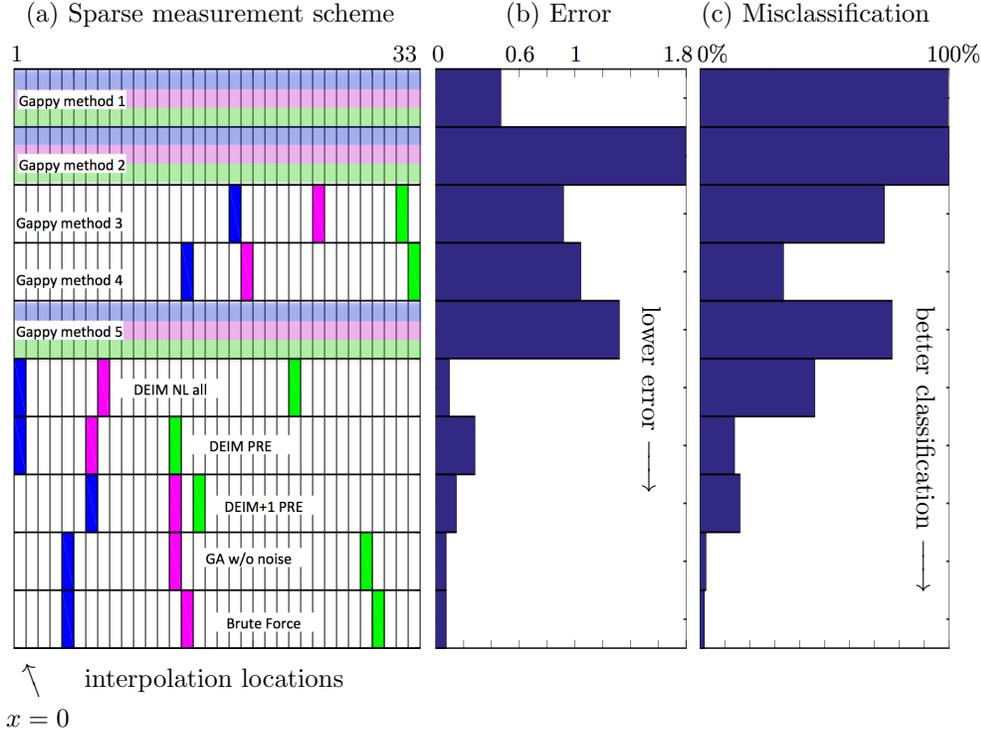
\begin{figure}[t]
\centering
\begin{overpic}[width=1.\textwidth]{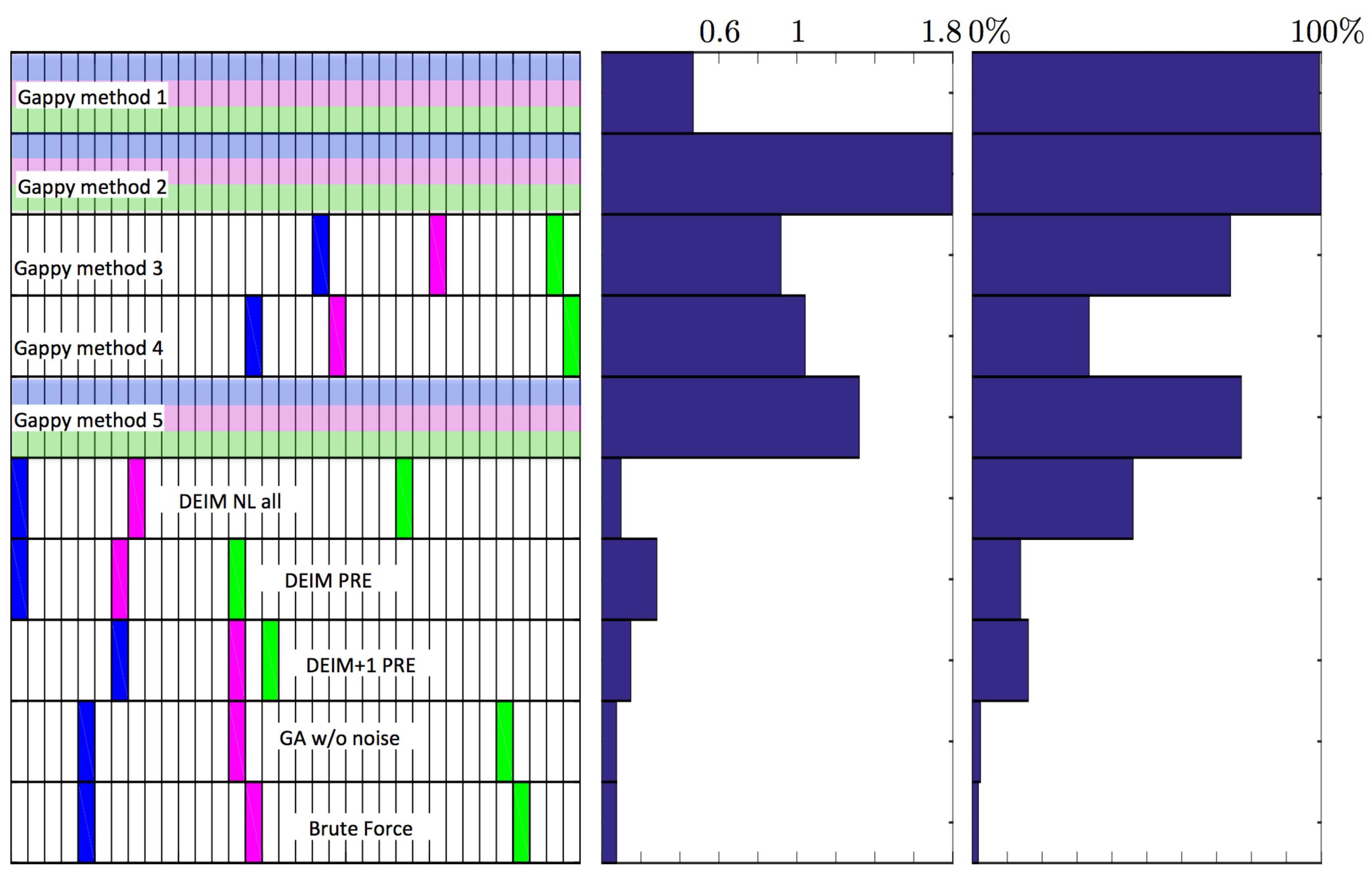}
\put(8,-3){interpolation locations}
\put(0,-7){$x=0$}
\put(1,-4){\rotatebox{110}{$\xrightarrow{\hspace*{0.3cm}}$}}
\put(43.5,61){\small{0}}
\put(0.5,61){\small{1}}
\put(39.5,61){\small{33}}
%\put(44,-2){(b)}
%\put(72,-2){(c)}
\put(65,36){\rotatebox{-90}{lower error $\xrightarrow{\hspace*{0.5cm}}$}}
\put(93,36){\rotatebox{-90}{better classification $\xrightarrow{\hspace*{0.5cm}}$}}
\put(51,65){(b) Error}
\put(71,65){(c) Misclassification}
\put(2,65){(a) Sparse measurement scheme}
\end{overpic}
\vspace*{.1in}
\caption{(a) Comparison of various sparse sampling strategies, including various gappy POD and DEIM methods, against 
the optimal sparse sampling solution determined by exhaustive search.  The first 33 indices of the
discretized PDE are shown where the first index corresponds to $x=0$.  The first, second and third
measurement locations are denoted by the blue, magenta and green bars respectively.   
The color bars that span the index locations represents random measurement locations which can be potentially at any of the indices.
The accompanying (b) error and (c) misclassification scores are also given.   In contrast to many of the other methods, the genetic algorithm proposed produces results that are almost identical to the true optimal solution, making it a viable method for online
enhancement of ROMs.  The various sparse selection methods are as follows:  (i) Gappy method 1:  random selection of three
indices from interpolation range 1-33 (where the histograms in Figs.~\ref{fig:cqg}-\ref{fig:gl2} suggest the measurements should occur), (ii) 
Gappy method 2: random selection from all possible interpolation points on the domain, (iii) Gappy method 3:  condition number minimization routine for three interpolation points~\cite{gap2}, (iv) Gappy method 4:  same as Gappy method 3 but with ten interpolation points (i.e. it is now full rank), (v) Gappy method 5:  selection of interpolation points from maxima and minima of POD modes~\cite{karni}, (vi) DEIM NL all:  DEIM algorithm applied jointly to all the nonlinear terms of all dynamical regimes, (vii) DEIM PRE:  algorithm developed in \cite{sargsyan2015}, (viii) DEIM+1 PRE:  use the algorithm in DEIM PRE but discard the first DEIM point
and select from the 2nd-4th DEIM points, (ix) GA:  genetic algorithm advocated here, (x) Brute force:  optimal solution from exhaustive search }
%
%Right - Gappy method 1 - randomly choose three locations from 0:33 (this is the range of locations from PRE when running DEIM for different regimes). Gappy Method 2 - randomly select three locations from 0:512;
%Gappy Method 3 - wilconx 3;
%Gappy Method 4 -  wilconx 10; 
%Gappy Method 5 - Karni; 
%DEIM NL all - DEIM applied to nonlinear terms from  all regimes together; 
%DEIM PRE; 
%DEIM + 1 PRE - run algorithm from PRE paper with 4 sensor locations then discard the first location; 
%GA  - genetic algorithm; 
%Middle - relative reconstruction error; 
%Right - misclassification }
 \label{fig:gl3}
\end{figure}

To finish our analysis, we compare the DEIM architecture against some classic gappy POD and DEIM methods.  Figure~\ref{fig:gl3}
gives an algorithmic comparison of the interpolation point selection of various techniques against the proposed method
and the ground truth optimal solution obtained by exhaustive search.   Both the reconstruction error and classification accuracy
are important in selecting the interpolation indices, and both are represented in panels (b) and (c) of Fig.~\ref{fig:gl3}.
Importantly, the method proposed here, which starts with the DEIM+1 points and does a quick genetic algorithm search
produces nearly results that are almost equivalent to the exhaustive search. This is an impressive result given the efficiency 
of the genetic search and online computation possibilities.  And even if one is not interested in executing the genetic
search, the DEIM+1 points used for ${\bf P}$ provide nearly double the performance (in terms of accuracy) versus DEIM.

\subsection{Flow Around a Cylinder}

The previous example provides an excellent proof-of-concept given that we could compute a ground truth optimal
solution.  The results suggest that we should start with the DEIM measurement matrix ${\bf P}$ and execute
the genetic algorithm from there.  We apply this method on the classic problem of flow around a cylinder.  
This problem is also well understood and has already been the subject of
studies concerning sparse spatial measurements~\cite{Bright:2013,gappy,lionel,lionel2,bing}.  
Specifically, it is known that for low to moderate Reynolds numbers, the dynamics is spatially
low-dimensional and POD approaches have been successful in quantifying the
dynamics~\cite{gappy,pod1,pod2,pod3,pod4}. 
The Reynolds number, $Re$, plays the role of the bifurcation parameter $\mu$ in (\ref{eq:complex}), i.e. it is a parametrized dynamical system.

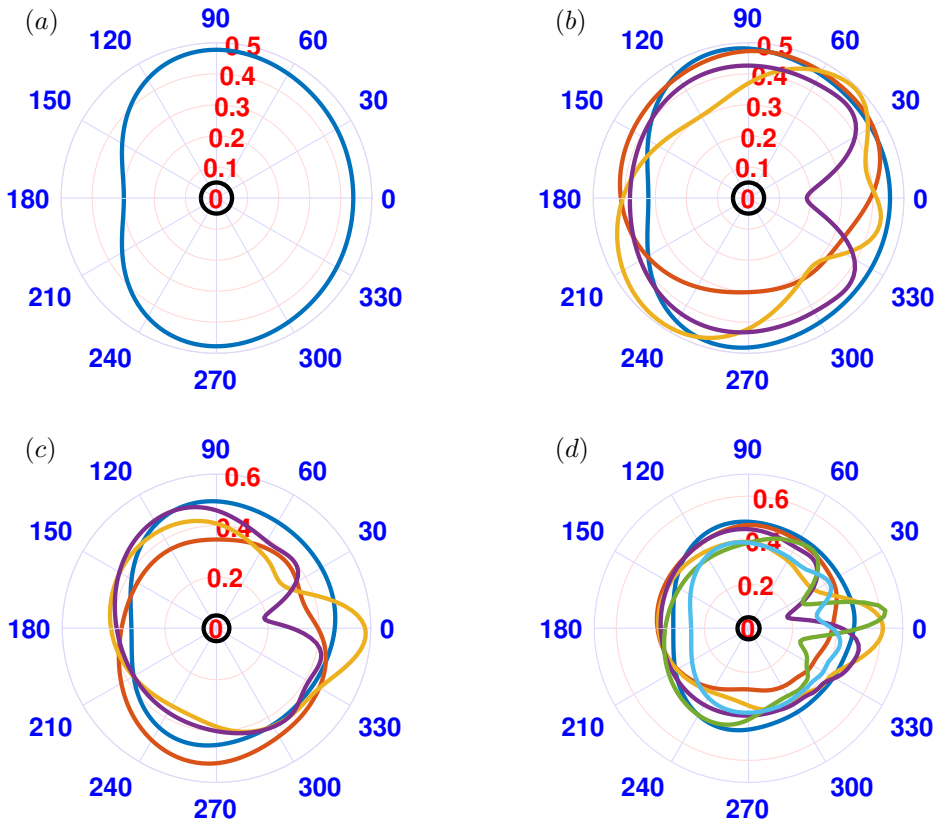
\begin{figure}[t]
%\hspace*{-.6in}
%\begin{overpic}[width=1.25\textwidth]{cylinder_modes04.eps}
%\hspace*{1.5in}
\begin{overpic}[width=0.95\textwidth]{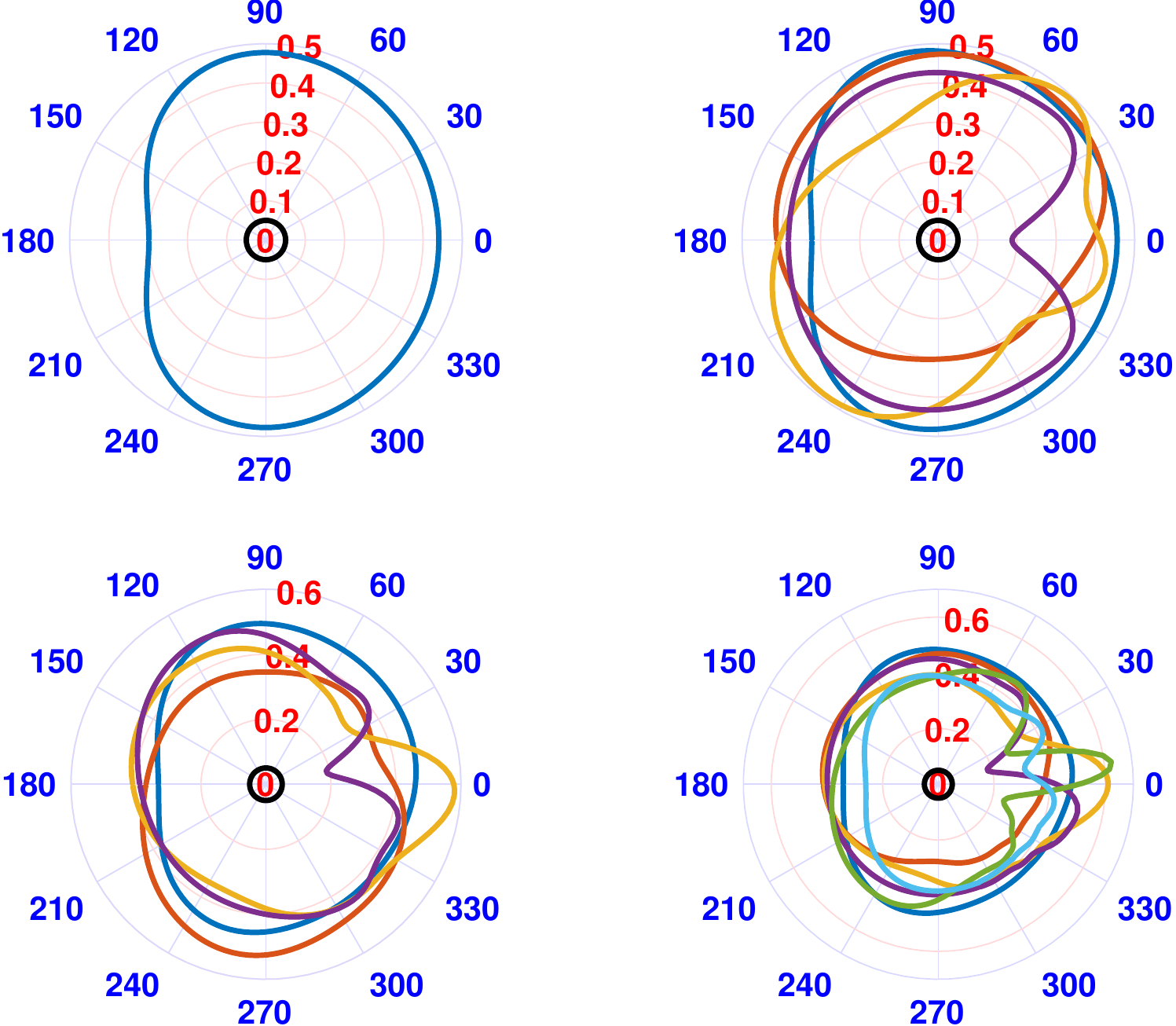}
\put(2,85){$(a)$}
\put(59,85){$(b)$}
\put(2,39){$(c)$}
\put(59,39){$(d)$}
\end{overpic}
\vspace*{-.0in}
\caption{Dominant POD modes for Reynolds numbers (a) 40, (b) 150, (c) 300 and (d) 1000.  The
POD modes are plotted in cylindrical coordinates to illustrate the pressure field generated on the cylinder.  
For low Reynolds numbers ($Re=40$), a single dominant POD mode exists.  As the Reynolds number
increases, more exotic mode structures, which are asymmetric, are manifested.  The blue labels are
the degrees around the cylinder while the red labels are the amplitudes of the POD pressure modes.  The
first mode is in blue, followed by red, gold, purple, green and cyan.}
 \label{fig:cyl1}
\end{figure}

The data we consider comes from numerical simulations of the incompressible
Navier-Stokes equation: 
\begin{subeqnarray}
&&\frac{\partial u}{\partial t}+u\cdot\nabla u+\nabla p-\frac{1}{Re}\nabla^{2}u=0\\
&&\nabla\cdot u=0
\label{eq:incompresNS}
\end{subeqnarray}
where $u\left(x,y,t\right)\in \mathbb{R}^{2}$ represents the 2D velocity, and
$p\left(x,y,t\right)\in \mathbb{R}^2$ the corresponding pressure field. The boundary condition are
as follows: (i) Constant flow of $u=\left(1,0\right)^{T}$ at $x=-15$, i.e., the
entry of the channel, (ii) Constant pressure of $p=0$ at $x=25$, i.e., the end of the channel, and (iii) 
Neumann boundary conditions, i.e. $\frac{\partial u}{\partial\mathbf{n}}=0$
on the boundary of the channel and the cylinder (centered at $(x,y)=(0,0)$ and of radius unity). 

We consider the fluid flow for Reynolds number $Re=40,150,300,1000$ and perform an SVD
on the data matrix in order to extract POD modes.   The rapid decay of singular values allows us
to use a small number of POD modes to describe the fluid flow and build local ROMs.  The POD modes retained for
each Reynolds number is shown in Fig.~\ref{fig:cyl1}.  These modes are projected on cylindrical
coordinates to better demonstrate the structure of the pressure field generated on the cylinder.

The POD modes can be used to construct a DEIM interpolation matrix ${\bf P}$ illustrated in Fig.~\ref{fig:ga}.
The DEIM interpolation points already provide a good set of interpolation points for classification
and reconstruction of the solution.  
However, the genetic algorithm advocated in this work can be used to adjust
the interpolation points and achieve both better classification performance and improved reconstructions.
In the cubic-quintic GL equation, the error was improved by nearly an order of magnitude over the standard
DEIM approach.  For the flow around the cylinder, the error is also improved from the DEIM algorithm, quickly
reducing the error with $M=2$ generations and converging to the nearly optimal interpolation points within $M=10$ generations.
Given the limited number of interpolation points, the genetic search can be computed in an online
fashion even for this two-dimensional fluid flow problem.

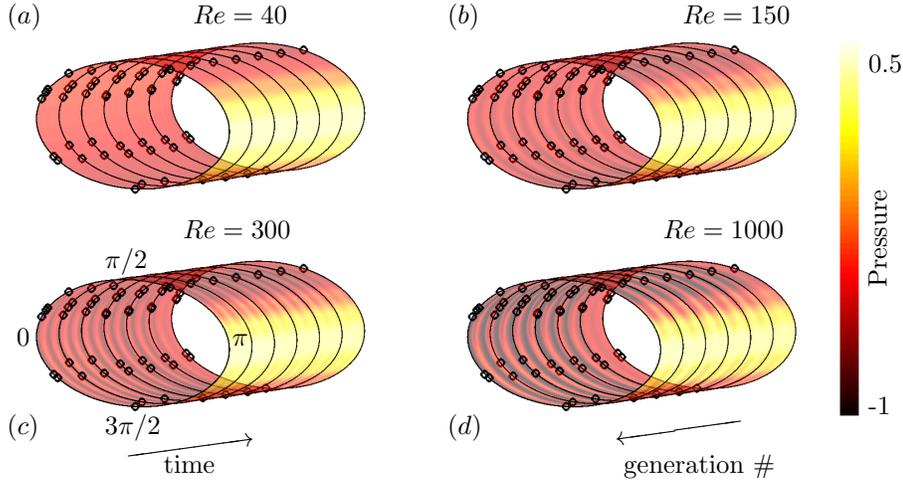
\begin{figure}[t]
%\begin{overpic}[width=0.5\textwidth]{cylinder_9sens.eps}
%\end{overpic}
%\hspace*{-.2in}
\begin{overpic}[width=1\textwidth]{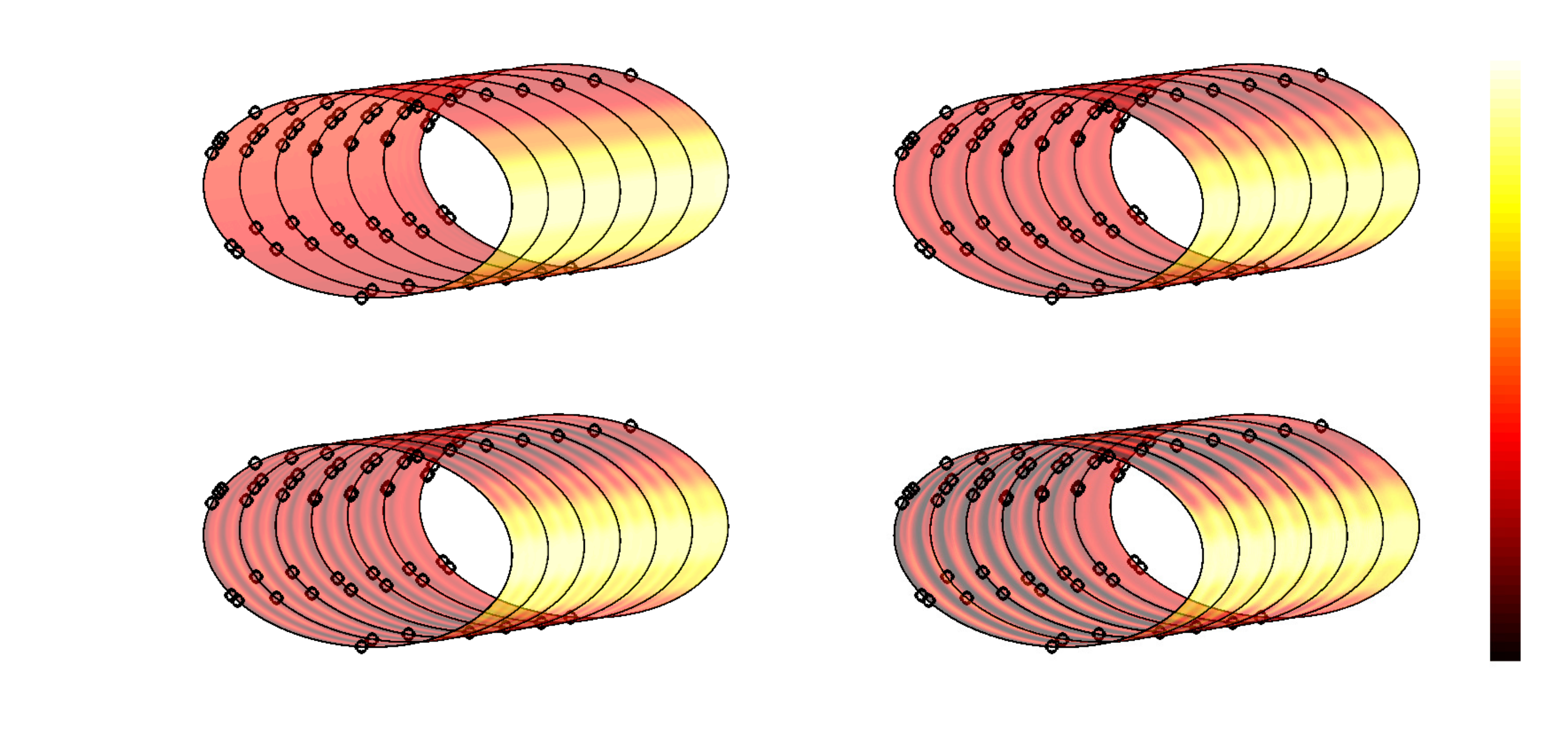}
%\put(50,55){$\beta_3$}
\put(10,45){$(a)$}
\put(55,45){$(b)$}
\put(10,3){$(c)$}
\put(55,3){$(d)$}
\put(29,45){$Re=40$}
\put(28,23){$Re=300$}
\put(79,45){$Re=150$}
\put(77,23){$Re=1000$}
\put(11,12){0}
\put(20,20){$\pi/2$}
\put(33,12){$\pi$}
\put(20,3){$3\pi/2$}
\put(26,-1){time}
%\put(15,35){0}
%\put(28,37.5){$T$}
\put(22,0){\rotatebox{8}{$\xrightarrow{\hspace*{1.5cm}}$}}
\put(72,5){\rotatebox{-172}{$\xrightarrow{\hspace*{1.5cm}}$}}
\put(73,-1){generation \#}
\put(98,18){\rotatebox{90}{Pressure}}
\put(98,5){-1}
\put(98,40){0.5}
\end{overpic}
\caption{The heat map on the cylinder shows the dominant, low-dimensional pattern of activity that is used for generating POD modes
for Reynolds number (a) 40, (b) 150, (c) 300 and (d) 1000.  Overlaid on the heat map are the best sensor/interpolation locations (circles) at each generation of the genetic algorithm scheme for $m=10$ interpolation points over $M=7$ generations of the search.
The interpolation locations for each generation in the genetic algorithm start at the back and move forward.}
 \label{fig:cyl3}
\end{figure}

\begin{figure}[t]
%\begin{overpic}
%[width=0.55\textwidth]{cylinder_GA9.eps}
%\end{overpic}
%\hspace*{-.3in}
\begin{overpic}
[width=0.55\textwidth]{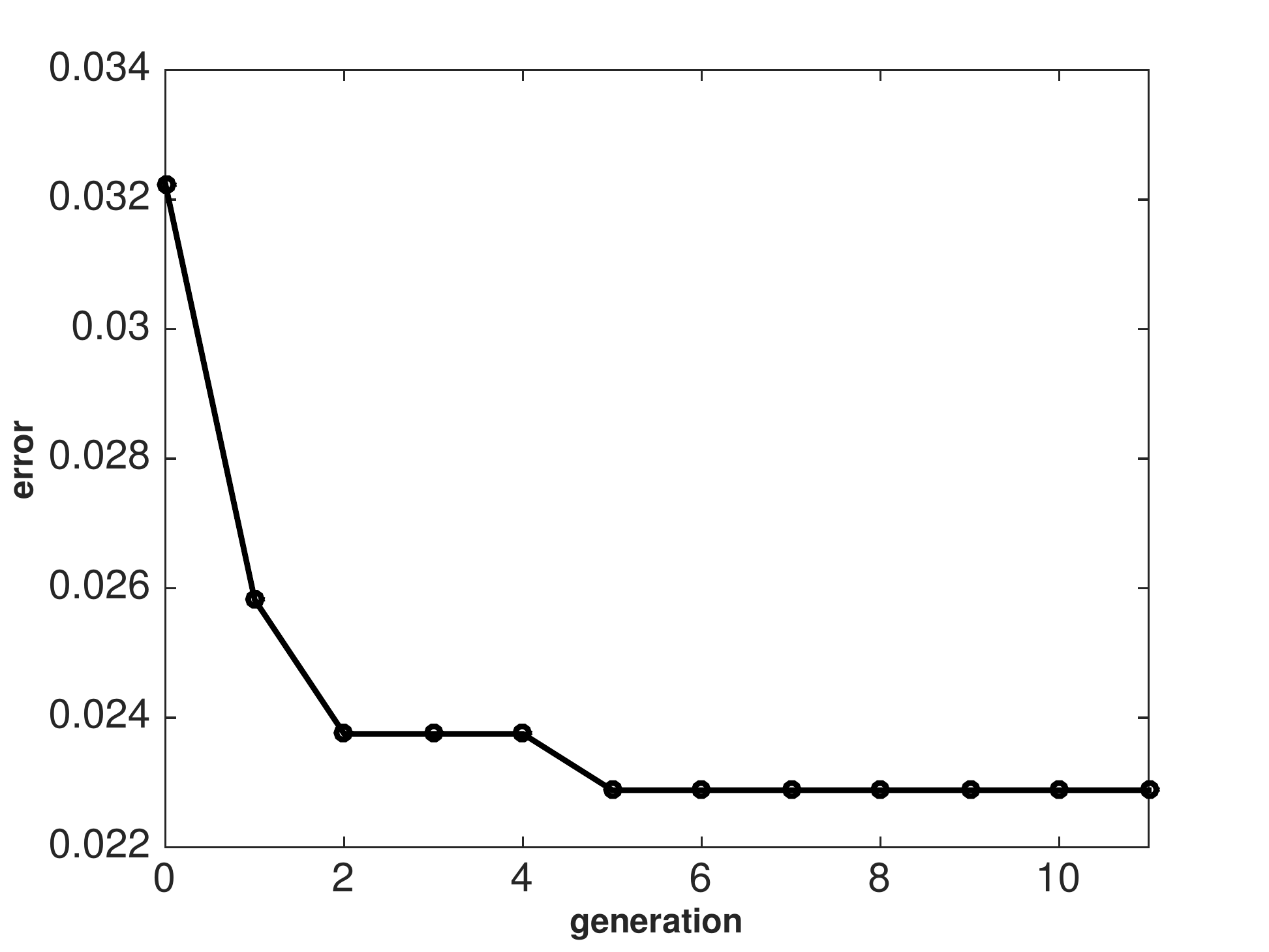}
\put(16,63){$(a)$}
\end{overpic}
\hspace*{-.3in}
\begin{overpic}
[width=0.55\textwidth]{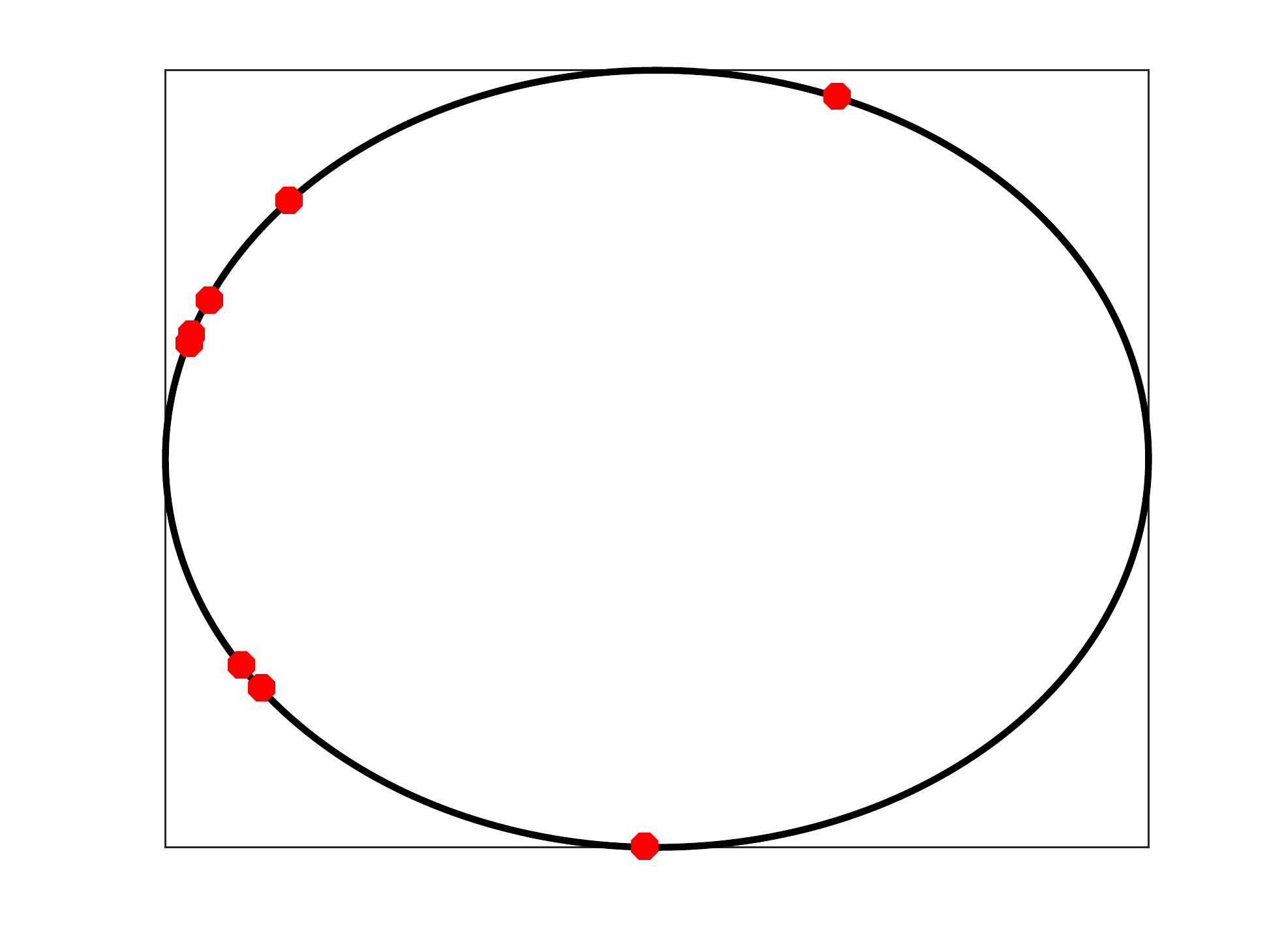}
%\put(50,55){$\beta_3$}
\put(16,63){$(b)$}
\end{overpic}
\vspace{-.25in}
\caption{ (a) Error reduction and convergence of the genetic algorithm from an initial DEIM+1 interpolation matrix ${\bf P}$.  Significant
error reduction is accomplished in a single generation.  A total of $M=10$ generations are shown, illustrating the convergence to the 
nearly optimal solution.  (b) The final position of the interpolation points (red) is shown for 10 interpolation points on the cylinder.}
 \label{fig:cyl4}
\end{figure}

Figure~\ref{fig:cyl3} is a composite figure showing the pressure field evolution in time along a the cylinder.
The heat map shows the dominant, low-dimensional pattern of activity that is used for generating POD modes.
Overlaid on this are the best sensor/interpolation locations at each generation of the genetic algorithm scheme
for 10 interpolation points over 7 generations of the search.  Note the placement of the interpolation points 
around the cylinder.  Specifically, as the number of generations increases, the interpolation points move to
better sampling positions, reducing the error in the ROM.  The convergence of the error across 10 generations
of the algorithm is shown in Fig.~\ref{fig:cyl4} along with the final placement of the interpolation points.  The near
optimal interpolation points are not trivially found.  Overall, the DEIM architecture with genetic algorithm search
reduces the error by anywhere between a factor of two and an order of magnitude, making it attractive for
online error reduction and ROM performance enhancement.

%\begin{figure}[t]
%%\centering
%\begin{overpic}
%[width=0.55\textwidth]{cylinder_sensors9.eps}
%\put(50,55){$\beta_3$}
%\end{overpic}
%\hspace*{-.3in}
%\begin{overpic}
%[width=0.55\textwidth]{cylinder_sensors10.eps}
%\put(50,55){$\beta_3$}
%\end{overpic}
%\caption{Left: For each regime we used 9 sensor locations then built the histogram of those locations. There were 7 (out of 9) locations that got more than one vote and we picked those. Right: The same with 10 sensors. We used 9 and 10 sensors since only those gave correct classification(without noise) results.}
% \label{figure10}
%\end{figure}

\section{Conclusions}

ROMs are enabled by two critical steps:  (i) the construction of a low-rank subspace where the dynamics can be accurately projected, and (ii) a sparse sampling method that allows for an interpolation-based projection to provide a nearly $\ell_2$ optimal subspace approximation to the nonlinear term without the expense of orthogonal projection.   Innovations to improve these two
aims can improve the outlook of scientific computing methods for modern, high-dimensional simulations that are rapidly
approaching exascale levels.  
These methods also hold promise for real-time control of complex systems, such as turbulence~\cite{Brunton2015amr}.  
This work has focused on improving the sparse sampling method commonly used in
the literature.  In partnership with the DEIM infrastructure, a genetic algorithm was demonstrated to determine 
nearly optimal sampling locations for producing a
subspace approximation of the nonlinear term without the expense of orthogonal projection.  The algorithm can
be executed in a potentially online manner, improving the error by up to an order-of-magnitude in the examples 
demonstrated here.   In our complex cubic-quintic Ginzburg-Landau equation example,  for a fixed number of interpolation points $m$, the first $m+1$ DEIM interpolation points
are computed and the first point is discarded.  This DEIM+1 sampling matrix alone can reduce the error by a factor of two before starting the genetic algorithm search.

In general, genetic algorithms are not ideal for optimization since they rarely have guarantees on convergence and have
many potential pitfalls.  In our case, the DEIM starting point for the interpolation point selection algorithm
is already close to the true optimum.  Thus the genetic algorithm is not searching blindly in a high-dimensional
fitness space.  Rather, the algorithm aims to simply make small adjustments and refinements to the sampling matrix in order
to maximize the performance of the nonlinear interpolation approximation.   In this scenario, many of the
commonly observed genetic algorithm failures are of little concern.   The method is shown to reduce the
error by a substantial amount within only one or two generations, thus making it attractive for implementation in 
real large-scale simulations where accuracy of the solution may have significant impact on the total computational
cost.  In comparison to many other sparse sampling strategies used in the literature, it out performs them by a
significant amount both in terms of accuracy and ability to classify the dynamical regime.  Indeed, the algorithm
refines the sampling matrix ${\bf P}$ to be nearly optimal.

\section*{Acknowledgements}
J. N. Kutz would like to acknowledge support from the Air Force Office of Scientific Research (FA9550-15-1-0385).

\end{document}